\documentclass[a4paper,twoside,10pt]{amsart}

\usepackage{bbm}
\usepackage{psfrag}
\usepackage{graphicx}
\usepackage{amssymb}
\usepackage{amsmath}
\usepackage{amscd}
\usepackage{a4wide}
\usepackage[colorlinks,linkcolor=blue,citecolor=blue,urlcolor=blue,pdfpagemode=UseOutlines,bookmarks]{hyperref}

\newcommand{\complex}[1]{\mathsf{#1}}

\newcommand{\Boxcomplex}[1]{\mathsf{\complex B}_{0}(#1)}

\newcommand{\geom}[1]{\| #1 \|}

\newcommand{\opcd}[3]{\mathrm{cd}_{#2}^{#1}#3}

\newcommand{\zind}[2]{\mathrm{ind}_{\Z_{#1}}(#2)}

\newcommand{\RS}[3]{#3_{#2}^{#1}}

\newcommand{\sd}[1]{\mathrm{sd}(#1)}

\newcommand{\id}{\mathop{\rm Id}}

\newcommand{\Dim}[1]{\mathop{\rm dim}(#1)}

\newcommand{\im}[1]{{\mathop{\rm Im}#1}}

\newcommand{\vertices}[1]{\mathrm{V}(#1)}

\newcommand{\edges}[1]{\mathrm{E}(#1)}

\newcommand{\chr}[1]{\chi(#1)}

\newcommand{\KG}[3]{\mathrm {KG}_{#2}^{#1}#3}
 
\newcommand{\kg}[3]{\mathrm {kg}_{#2}^{#1}#3}

\newcommand{\set}[2]{\left\{#1\vphantom{#2}\right.\;\left|\;\vphantom{#1}#2\right\}}

\newcommand{\Z}{\mathbbm{Z}}
\newcommand{\R}{\mathbbm{R}}

\newcommand{\s}{\mathcal{S}}
\newcommand{\T}{\mathcal{T}}

\newtheorem{theorem}{Theorem}

\begin{document}

\title{On generalised Kneser colourings}

\author{Carsten Lange}
\address{TU Berlin\\
        Institut f\"ur Mathematik, MA 6-2\\
	Stra\ss e des 17. Juni 136\\
	10623 Berlin}
\email{lange@math.tu-berlin.de}
\thanks{Supported by the DFG Sonderforschungsbereich 288
                          ``Differentialgeometrie und Quantenphysik'' in Berlin.}

\date{December 3, 2003}

\begin{abstract}
There are two possible definitions of the ``$s$-disjoint $r$-uniform Kneser hypergraph'' of a set system $\T$: The hyperedges
are either $r$-sets or $r$-multisets. We point out that Ziegler's (combinatorial) lower bound on the chromatic number of 
an $s$-disjoint $r$-uniform Kneser hypergraph only holds if we consider $r$-multisets as hyperedges. We give a new proof 
of his result and show by example that a similar result does not hold if one considers $r$-sets as hyperedges.  

In case of $r$-sets as hyperedges and $s \geq 2$ the only known lower bounds are obtained from topological invariants of associated
simplicial complexes if $r$ is a prime or the power of prime. This is also true for arbitrary $r$-uniform hypergraphs with $r$-sets
or $r$-multisets as hyperedges as long as $r$ is a power of a prime.
\end{abstract}

\maketitle

\section{Introduction}

During the last 25 years, topological methods have been successfully applied to obtain lower bounds of the chromatic number
of graphs and hypergraphs. It all took off with Lov{\'a}sz' seminal paper~\cite{L} were he proved that the connectivity of the 
neighbourhood complex of a graph~$G$ can be used to establish a lower bound of its chromatic number~$\chr G$ and that this bound 
can be used to prove Kneser's conjecture. The following diagram indicates how Lov{\'a}sz' result for Kneser graphs has been 
generalised in the subsequent years to Ziegler's result for $s$-disjoint $r$-uniform Kneser hypergraphs~$\KG{r}{s}{\T}$ associated 
to a family~$\T \subseteq 2^{[n]}$ of subsets of an $n$-set $[n]$. 
{\footnotesize
\[
  \begin{CD}
     \begin{matrix}
        \text{Lov{\'a}sz~\cite{L}}\\
	\opcd{2}{1}{\binom{[n]}{k}} \leq \chr {\KG{2}{1}{\binom{[n]}{k}}}
     \end{matrix}
     @>>> 
     \begin{matrix}
        \text{Alon, Frankl, Lov{\'a}sz~\cite{AFL}}\\
	\opcd{r}{1}{\binom{[n]}{k}} \leq (r-1) \cdot \chr {\KG{r}{1}{\binom{[n]}{k}}}
     \end{matrix}
     @>>> 
     \begin{matrix}
        \text{Sarkaria~\cite{S90}}\\
     	\opcd{r}{s}{\binom{[n]}{k}} \leq (r-1) \cdot \chr {\KG{r}{s}{\binom{[n]}{k}}}
     \end{matrix}\\
     @VVV @VVV @VVV \\
     \begin{matrix}
        \text{Dol'nikov~\cite{D81}}\\
	\opcd{2}{1}{\T} \leq \chr {\KG{2}{1}{\T}}
     \end{matrix}
     @>>> 
     \begin{matrix}
        \text{K{\v r}{\'i}{\v z}~\cite{K92,K00}}\\
	\opcd{r}{1}{\T} \leq (r-1) \cdot \chr {\KG{r}{1}{\T}}
     \end{matrix}
     @>>> 
     \begin{matrix}
        \text{Ziegler~\cite{Z}}\\
     	\opcd{r}{s}{\T} \leq (r-1) \cdot \chr {\KG{r}{s}{\T}}
     \end{matrix}
  \end{CD}
\]
}We defer precise definitions of the $s$-disjoint $r$-colourability defect $\opcd{r}{s}{\T}$ and of $\KG{r}{s}{\T}$ to
Section~\ref{sec_prelims}. For $s > 1$, the Kneser hypergraph $\KG{r}{s}{\T}$ just mentioned is not a hypergraph in the sense of 
Berge~\cite{Berge}. It is a hypergraph \emph{with multiplicities}, that is, a hyperedge may contain multiple copies of 
nodes but must contain at least two different nodes. Usually a hyperedge of a hypergraph consists of distinct nodes, 
i.e.\ the hypergraph is \emph{multiplicity-free}. In Section~\ref{sec_prelims} we also define multiplicity-free Kneser 
hypergraphs $\kg{r}{s}{\T}$ (a hypergraph in the sense of Berge) and point out that the proof of Theorem~5.1 of Ziegler~\cite{Z} 
does not work in the generality stated there. Section~\ref{sec_examples_and_counterexamples} is devoted to examples. We show 
that neither Sarkaria's nor Ziegler's result holds if we replace $\KG{r}{s}{\T}$ by $\kg{r}{s}{\T}$.

In Section~\ref{sec_prime} we extend a result from Alon, Frankl, and Lov{\' a}sz~\cite{AFL} for multiplicity-free $r$-uniform 
hypergraphs : We give a topolgical lower bound of the chromatic number of an $r$-uniform hypergraph with or without multiplicities 
if $r$ is the power of a prime. A review of known facts about free and fixed-point free group actions and simplicial complexes is 
given in Section~\ref{sec_group_actions}. We also describe a variety of simplicial complexes that we need later, such as a box 
complex associated to an arbitrary $r$-uniform hypergraph $\s$. This complex is different from (but closely related to) the complex 
considered by Alon, Frankl, and Lov{\'a}sz. In an unpublished preprint from 1987, \"Ozaydin extended their result to the case where 
$r$ is the power of a prime,~\cite{O87}.

In Section~\ref{sec_loops} we give a new proof of Ziegler's result for $r$-uniform $s$-disjoint Kneser hypergraphs with multiplicities.
The main tool used is ``Sarkaria's inequality''. The proof is inspired by Matou{\v s}ek's proof~\cite{M02} of the result of 
K{\v r}{\' i}{\v z}~\cite{K92,K00}. In fact we prove a bit more: The colourability defect of a set system $\T$ is not the 
only lower bound for the chromatic number of the $r$-uniform Kneser hypergraph associated to $\T$. If $r$ is prime we squeeze between 
these two numbers the index of an associated simplicial complex. Theoretically, this might yield better lower bounds but has not been 
investigated so far.

For arbitrary $r$-uniform hypergraphs with or without multiplicities, no lower bound except the topological one stated in 
Theorem~\ref{thm_primepower} is known. The problem one is confronted with is the same that one faces in proving the topological 
Tverberg theorem. Only partial results in case of prime-powers are achieved so far. Moreover, it is not known whether every 
$r$-uniform hypergraph (multiplicity-free or not) can be realised as an $r$-uniform $s$-disjoint Kneser hypergraph for some $s$ 
and an appropriate set system. A result of this type is known for graphs: Every graph is a Kneser graph,~\cite{MZ}.

\section{Preliminaries}
\label{sec_prelims}

\noindent
\textbf{$\mathbf{s}$-disjoint sets.}\ \
For integers $r > s \geq 1$ we say that the subsets $S_{1}, {\ldots} , S_{r}$ of~$[n]$ 
are $s$\emph{-disjoint} if each element of~$[n]$ occurs in at most~$s$ of the sets~$S_{i}$, or equivalently, if the 
intersection of any choice of~$s+1$ sets is empty. This formulation is the reason that this concept is called $s$-wise 
disjoint by Sarkaria~\cite{S90}. We emphazise that $S_{i}= S_{j}$ may occur for $i \neq j$.

\noindent
\textbf{$\mathbf{s}$-disjoint $\mathbf{r}$-colourability defect.}\ \
The $s$\emph{-disjoint} $r$\emph{-colourability defect} $\opcd{r}{s}{\s}$ of $\s \subseteq 2^{[n]}\setminus \emptyset$
is the number of elements one has to remove from the multiset~$[n]^{s}$ such that the remaining multiset can be covered by an
$s$-disjoint $r$-tuple of sets such that none of the sets contains an element from~$\s$. This number can be computed by evaluating
\[
  \opcd{r}{s}{\s} = n \cdot s - \max \set{\sum_{j=1}^{r}|R_{j}|}
                                         {\begin{matrix}
                                             R_{1}, {\ldots} , R_{r} \subseteq [n]\ s\text{-disjoint}\\
					     \text{and } S \not \in R_{j} \text{ for all } S \in \s \text{ and all }j
					  \end{matrix}}.
\]

\noindent
\textbf{$\mathbf r$-multisubsets of ${\mathbf {[n]}}$.}\ \
The collection $x_{1}, {\ldots} , x_{r}$ of elements of $[n]$ is called an $r$\emph{-multisubset} of~$[n]$ if the set 
$\{ x_{1},{\ldots} ,x_{r} \}$ has cardinality at least two. We denote an $r$-multiset by $\{\!\{ x_{1},{\ldots} ,x_{r} \}\!\}$.

\noindent
\textbf{$\mathbf r$-uniform hypergraphs with or without multiplicities.}\ \
Consider $\s \subseteq 2^{[n]}$ such that $\bigcup_{S \in \s}S = [n]$. An $r$-uniform hypergraph $\s$ without multiplicities coincides
with Berge's definition of an $r$-uniform hypergraph, \cite{Berge}: The vertices are $[n]$ and the hyperedges are the 
$r$-subsets~$\s$ of~$[n]$.  
Let~$\s^{\prime}$ be a set of $r$-multisubsets of~$[n]$ such that for every $i \in [n]$ exists an $S \in \s^{\prime}$ with $i \in S$. 
We call~$\s^{\prime}$ an $r$\emph{-uniform hypergraph with multiplicities} and say that it has \emph{node set} 
$\vertices {\s^{\prime}} = [n]$ and \emph{hyperedge set} $\edges {\s^{\prime}} = \s^{\prime}$.

\noindent
\textbf{$\mathbf r$-uniform $\mathbf s$-disjoint Kneser hypergraphs.}\ \
For a set $\T = \{ T_{1},{\ldots} ,T_{m} \}$ of subsets of~$[n]$, we consider the $r$\emph{-uniform} $s$\emph{-disjoint}
\emph{Kneser hypergraph} $\KG{r}{s}{\T}$ \emph{with multiplicities} on the node set $\vertices {\KG{r}{s}{\T}} = [m] = [| \T |]$ 
with hyperedges
\[
  \edges {\KG{r}{s}{\T}} := \set{ \ \{\!\{ k_{1},{\ldots}, k_{r} \}\!\} }
                                { \begin{matrix}
                                     k_{i} \in [m],\ \{\!\{ k_{1},{\ldots}, k_{r} \}\!\} \text{ is an $r$-multiset, }\\ 
   				     \text { and }T_{k_{1}},{\ldots}, T_{k_{r}} \text{ are $s$-disjoint }
				  \end{matrix}
				}.
\]
The $r$\emph{-uniform} $s$\emph{-disjoint Kneser hypergraph} $\kg{r}{s}{\T}$ \emph{without multiplicities} has the same node set
$\vertices {\kg{r}{s}{\T}} = [m] =[ | \T | ]$ but a different hyperedge set:
\[
  \edges {\kg{r}{s}{\T}} := \set{ \ \{ k_{1},{\ldots}, k_{r} \} \subseteq [m]}
                                { \begin{matrix}
                                    \{ k_{1},{\ldots}, k_{r} \} \text{ is an $r$-set }\\
				    \text{and } T_{k_{1}},{\ldots}, T_{k_{r}} \text{ are $s$-disjoint}
				  \end{matrix}
				}.
\]
In the special case $s = 1$ we have $\KG{r}{s}{\T} = \kg{r}{s}{\T}$ since an $r$-multiset with $r$ pairwise disjoint elements 
can be seen as an $r$-set. In particular for $r=2$ we have $s=1$ and both definitions specialise to a Kneser graph of 
$\T \subseteq 2^{[n]}$.

\noindent
\textbf{Colourings.}\ \
A colouring of an $r$-uniform hypergraph~$\s$ with~$m$ colours is a map $c: \vertices \s {\rightarrow} [m]$ that assigns to each 
node of~$\s$ a colour so that no hyperedge is monochromatic, that is, for $e \in \edges \s$ we have $| \set{c(x)}{x\in e} | \geq 2$. 
The \emph{chromatic number} $\chr \s$ is the smallest number~$m$ such that a colouring of~$\s$ with~$m$ colours exists. Every 
hyperedge of $\kg{r}{s}{\T}$ is a hyperedge of $\KG{r}{s}{\s}$, hence:
\[
  \chr {\kg{r}{s}{\T}} \leq \chr {\KG{r}{s}{\T}}.
\]

\begin{figure}[b]
   \psfrag{1}{$1$}
   \psfrag{2}{$2$}
   \psfrag{3}{$3$}
   \psfrag{4}{$4$}
   \psfrag{5}{$5$}
   \psfrag{A}{$A$}
   \psfrag{B}{$B$}
   \psfrag{C}{$C$}
   \psfrag{D}{$D$}
   \psfrag{E}{$E$}
   \psfrag{F}{$F$}
   (a) \includegraphics[width=3.5cm]{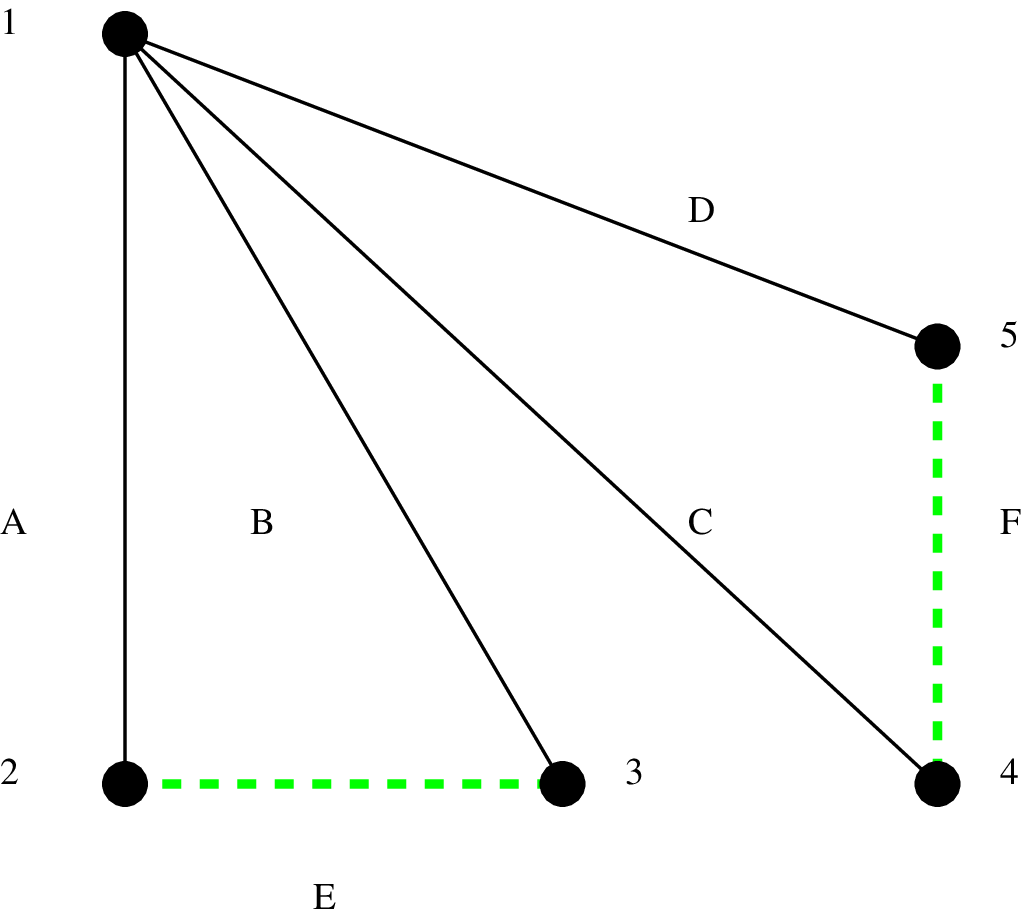}
   $\qquad\qquad$
   (b) \includegraphics[width=3.5cm]{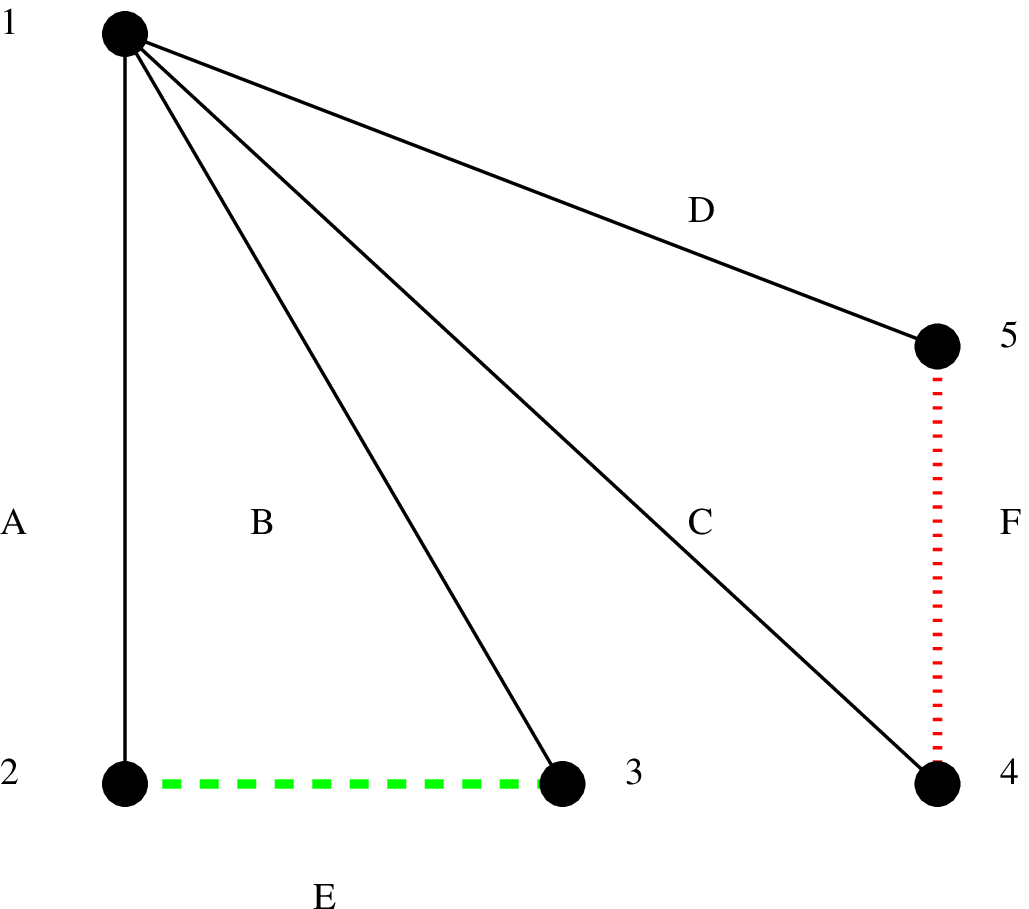}
   \caption{The edges $A$, $B$, $C$, $D$, $E$, and $F$ represent the nodes of $\kg{r}{s}{\T}$ and  $\KG{r}{s}{\T}$.
            In (a) the nodes of the hypergraph $\kg{r}{s}{\T}$ are coloured with two colours (best possible).
	    In (b) The nodes of the hypergraph $\KG{r}{s}{\T}$ are coloured with three colours (best possible).
	   }
   \label{illustration}
\end{figure}

\noindent
\textbf{A remark on Theorem~5.1 of~\cite{Z}.}\ \
Ziegler states his Theorem~5.1 for hypergraphs in the sense of Berge. The proof on page~679 of~\cite{Z}
yields the disired contradiction only if one assumes a colouring of a hypergraph with $r$-multisets as hyperedges. 
More precisely, the construction only guarantees that the~$p$ subsets $S_{1}, {\ldots}, S_{p}$ of~$[n]$ are $\boldmath s$-disjoint, 
they need not to be pairwise different. 

\section{Examples and Counterexamples}
\label{sec_examples_and_counterexamples}

\noindent
\textbf{Example 1.}\ \
We give an example to illustrate the two concepts of $s$-disjoint $r$-uniform Kneser hypergraphs and to see that the chromatic 
numbers $\chr {\kg{r}{s}{\s}}$ and $\chr {\KG{r}{s}{\s}}$ can be different. We restrict ourselves to the following small but 
interesting case: $r=3$, $s=2$, and $\T \subseteq \binom{[5]}{2}$. Let us consider  
\[
  \T := \{\  \{1,2\},\{1,3\},\{1,4\},\{1,5\},\{2,3\}, \{4,5\} \ \}.
\]
The $3$-uniform $2$-disjoint Kneser hypergraph $\kg{3}{2}{\T}$ without multiplicities has~$6$ nodes and any hyperedge consists 
either of two nodes of type $\{ 1,x \}$ plus $\{ 2,3\}$ or $\{ 4,5 \}$ or it consists of $\{ 2,3 \}$, $\{ 4,5 \}$ plus one node 
of type $\{ 1,x \}$. Colouring $\kg{r}{s}{\T}$ means therefore colouring the edges of the graph shown in Figure~\ref{illustration} 
such that no triple of edges that form a hyperedge is monochromatic. This can be done with~$2$ colours as indicated in 
Figure~\ref{illustration}: Colour each node that contains~$1$ with colour one, and colour $\{ 2,3 \}$ and $\{ 4,5 \}$ with the 
other colour. Among others, the $s$-disjoint $r$-uniform Kneser hypergraph $\KG{r}{s}{\T}$ with multiplicities has the following 
additional hyperedges:
\[
   \{\!\{\ \{2,3\}, \{2,3\}, \{4,5\} \ \}\!\}\qquad\text{ and }\qquad  \{\!\{\ \{2,3\}, \{4,5\}, \{4,5\} \ \}\!\}.
\]
These hyperedges force us to use at least three colours; for a colouring see also Figure~\ref{illustration}.
The sets $R_{1}=\{ 2,4 \}$, $R_{2}=\{ 2,5 \}$ and $R_{3} = \{ 3,4 \}$ are $2$-disjoint, i.e. $\opcd{r}{s}{\T} \leq 4$.
From Ziegler's theorem or our Theorem~\ref{thm_loop_kneser}, we know that 
$\frac{\opcd{3}{2}{\T}}{3-1} \leq \chr {\KG{3}{2}{\T}}$. In this particular example we have 
\[
  \frac{\opcd{3}{2}{\T}}{3-1} \leq 2 = \chr {\kg{3}{2}{\T}} < \chr {\KG{3}{2}{\T}} = 3.
\]

\noindent
\textbf{Counterexample 1.}\ \
Example 1 can be modified easily to show that the colourability defect is not a lower bound for $\kg{r}{r-2}{\T}$ in general. 
For fixed $n \geq 5$, we consider the following set~$\T$ of subsets of $[n]$:
\[
  \T := \left\{ \ \{ 1,2 \}, {\ldots} , \{ 1,n \}, \{ 2,3 \}, \{ 4,5 \} \ \right\}.
\]
The $(r-2)$-disjoint $r$-uniform Kneser hypergraph $\kg{r}{r-2}{\T}$ is easily described. Every hyperedge contains
$(r-2)$ different elements of $\left\{ \ \{ 1,2 \}, {\ldots} , \{ 1,n \} \ \right\}$ plus $\{ 2,3 \}$ and $\{ 4,5 \}$, so 
$\chr {\kg{r}{r-2}{\T}} = 2$ if $n \geq r-1$. To compute the $(r-2)$-disjoint $r$-colourability defect of $\T$, we cover $[n]^{r-2}$
by $r$ sets so that no set contains an element of $\T$. Obviously, a set that contains~$1$ does not contain any other element. Let 
$r_{1}$ denote the number of sets that contain~$1$. The sets that do not contain~$1$ cannot contain~$2$ and~$3$ as well as~$4$ and~$5$ 
at the same time. There are $r_{2} = r - r_{1}$ such sets. Therefore, we have not covered $(r-2) - r_{1}$ copies of~$1$, 
$2(r-2) - r_{2}$ copies of~$2$ or~$3$, and  $2(r-2) - r_{2}$ copies of~$4$ or~$5$. In other words, at least 
$(r-2) - r_{1} + 2(r-2) - r_{2} + 2(r-2) - r_{2} = 3n -10$ elements are not covered or $\opcd{r}{s}{\T} \geq 3r - 10$.
For $r > 8$ this implies $\opcd{r}{r-2}{\T} > (r-1)\chr{\kg{r}{r-2}{\T}}$.

\noindent
\textbf{Counterexample 2.}\ \
We now show that the colourability defect is not even a lower bound for $\kg{r}{r-1}{\binom{[n]}{2}}$ in general. 
From Ziegler~\cite[Lemma 3.2]{Z}, we know $\opcd{r}{s}{\binom{[n]}{2}} = \max \{ ns - r(k-1), 0 \}$, hence we have 
to show that $(r-1)\chr {\kg{r}{r-1}{\binom{[n]}{2}}} < (r-1)n - r$. It suffices to colour $\kg{r}{s}{\binom{[n]}{2}}$
with $n-2$ colours. This can be done by a greedy colouring already known to Kneser in case of graphs: Assign colour
$i$ to $T \in \T$ if $i$ is the smallest element of $T$ and $i \leq n-3$. The elements not coloured yet are $\{ n-2, n-1 \}$,
$\{ n-2, n \}$, and $\{ n-1,n \}$; too few to form a hyperedge. We colour these elements by colour $n-2$. This is in general
not optimal, but suffices.

\section{Groups acting on simplicial complexes}
\label{sec_group_actions}
This section summarises standard definitions and known facts. For many examples and a detailed treatment we refer to 
Matou{\v s}ek's textbook,~\cite{M03}.

\noindent
\textbf{Simplicial complexes.}\ \
An \emph{abstract simplicial complex} $\complex K$ is a finite hereditary set system with vertex set $\vertices {\complex K}$. 
For sets $A_{1},{\ldots} ,A_{t}$ we define 
$A_{1} \uplus {\ldots}  \uplus A_{t} := \set{(a,1)}{a\in A_{1}} \cup {\ldots} \cup \set{(a,t)}{a \in A_{t}}$.
An important construction in the category of simplicial complexes is the {\em deleted join operation}. 
For a simplicial complex $\complex K$ and positive integers $r \geq s$ the $r$\emph{-fold} $s$\emph{-wise deleted join}
${\complex K}^{*r}_{\Delta(s)}$ is defined as
\[
   \complex K^{*r}_{\Delta(s)} := \set{\complex F_{1} \uplus {\ldots}  \uplus \complex F_{r}}
                                      { \complex F_{i} \in \complex K \text{ and } 
		             	        \complex F_{1},{\ldots},\complex F_{r} \text{ is $s$-wise disjoint }}.
\]
To avoid confusion, we point to the fact that $s$ indicates $s$-wise disjointness in the definition not $s$-disjointness.
Any abstract simplicial complex $\complex K$ can be realized as a topological space $\| \complex K \|$ in~$\R^{d}$ for some $d$.

\noindent
\textbf{Free $\Z_{r}$-spaces and $\Z_{r}$-index.}\ \
A \emph{free $\Z_r$-space} is a topological space $X$ together with a free $\Z_{r}$-action $\Phi$, i.e. for all $g, h \in \Z_{r}$
we have $\Phi(g) \circ \Phi(h) = \Phi(g+h)$, $\Phi (0) = \id$, and $\Phi(g)$ has no fixed point for $g \in \Z_{r} \setminus \{ 0 \}$.
A continuous map $f$ between $\Z_{r}$-spaces $(X,\Phi_{X})$ and $(Y,\Phi_{Y})$ is \emph{$\Z_{r}$-equivariant} (or a 
\emph{$\Z_{r}$-map} for simplicity) if $f$ commutes with the $\Z_{r}$-actions, i.e. $f \circ \Phi_{X} = \Phi_{Y} \circ f$. A 
simplicial complex $(\complex K,\Phi)$ is a \emph{free simplicial $\Z_{r}$-space} (or a \emph{free simplicial $\Z_{r}$-complex}) if 
$\Phi: \complex K \to \complex K$ is a simplicial map such that $\geom \Phi $ is a free $\Z_{r}$-action on $\geom{ \complex K }$. 
A \emph{simplicial $\Z_{r}$-equivariant} map~$f$ is a simplicial map between two simplicial $\Z_{r}$-spaces that commutes with 
the free $\Z_{r}$-actions. An important class of free $\Z_{r}$-spaces is called $E_{n}\Z_{r}$-space: A free $\Z_{r}$-space is an
$E_{n}\Z_{r}$-space if it is $n$-dimensional and $(n-1)$-connected. The $\Z_{r}$-index $\zind{r}{X}$ of $(X,\Phi)$ is the
smallest $n$ such that a $\Z_{r}$-map from $X$ to some $E_{n}\Z_{r}$-space exists. A generalised Borsuk--Ulam theorem, e.g. Dold's 
theorem~\cite{Dold} for free group actions provides the index for $E_{n}\Z_{r}$-spaces: There is no $\Z_{r}$-map from 
$E_{n}\Z_{r}$ to $E_{n-1}\Z_{r}$. Since we consider cyclic group actions most of the time, we tend to refer to a $\Z_{r}$-space $X$
without explicit reference to $\Phi$.

\noindent
\textbf{Sarkaria's inequality.}\ \
A useful inequality concerning the $\Z_{r}$-index of the join $\complex K * \complex L$ of two free simplicial $\Z_{r}$-complexes 
$\complex K$ and $\complex L$ is \emph{Sarkaria's inequality}, \cite{M02,M03}: 
\[
  \zind{r}{\complex K * \complex L} \leq \zind{r}{\complex K} + \zind{r}{\complex L} + 1.
\]

\noindent
\textbf{Examples: $\mathbf{\RS{r}{s}{\complex P}}$ and $\mathbf{\RS{r}{s}{\complex P}\s}$.}\ \
Let $1 \leq s < r$ where $r$ is a prime.
Consider the poset $\RS{r}{s}{P}$ of $s$-disjoint $r$-tuples $(S_{1},{\ldots} ,S_{r})$ of subsets of~$[n]$ with
$\bigcup_{i \in [r]} S_{i} \neq \emptyset$ ordered by componentwise inclusion, that is, 
$(S_{1},{\ldots} ,S_{r}) \leq (T_{1},{\ldots} ,T_{r})$ if $S_{i} \subseteq T_{i}$ for all $1 \leq i \leq r$. The order complex 
of $\RS{r}{s}{P}$ is denoted by $\RS{r}{s}{\complex P}$ and can be interpreted as the barycentric subdivision of the $r$-fold 
$(s+1)$-wise deleted join of an $(n-1)$-simplex $\sigma^{n-1}$, which is homotopy equivalent to a wedge of $(ns -1)$-dimensional 
spheres,~\cite{M03} or alternatively~\cite{S90}. Hence $\RS{r}{s}{\complex P}$ is an $E_{ns-1}\Z_{r}$-space. Another space we shall need later is the 
order complex~$\RS{r}{s}{\complex P}\s$ of the subposet~$\RS{r}{s}{P}\s$ of~$\RS{r}{s}{P}$ that consists only of those $r$-tuples 
$(S_{1},{\ldots} ,S_{r})$ with $\sum_{i=1}^{r}| S_{i}| \geq n s - \opcd{r}{s}{\s} +1$, where $\s \subseteq 2^{[n]}$. Such tuples 
have the useful property that there exists an $S \in \s$ such that $S \subseteq S_{i}$ for some $i$.
 
\noindent
\textbf{Fixed-point free actions.}\ \
A \emph{fixed-point free action} $\Phi$ of a group $G$ of order $r$ on a topological space $X$ means that no $x \in X$ is fixed by 
all $g \in G$. Obviously, a free group action is also fixed-point free. A standard example is $(\R^{m})^{r}$
with the diagonal $\set{(v,{\ldots} ,v)}{v \in \R^{m}}$ removed. We have a standard action of $G$ on $(\R^{m})^{r}$ 
by permuting the $r$ copies of $\R^{m}$. The action is fixed-point free on 
$(\R^{m})^{r} \setminus \set{(v,{\ldots} ,v)}{v \in \R^{m}}$ for all $r$. It is free if $r$ is a prime. The space 
$(\R^{m})^{r} \setminus \set{(v,{\ldots} ,v)}{v \in \R^{m}}$ is homotopy
equivalent to an $(m(r-1) - 1)$-sphere. 
\begin{theorem}[``Volovikov's theorem'',~\cite{V96}]$ $\\ \label{thm_Volovikov}
   \noindent
   Let $r = p^{t}$ be a power of a prime and consider a fixed-point free action of~$(\Z_{p})^{t}$ on~$X$ and~$Y$. Suppose that
   for all $i \leq \ell$ we have $\widetilde{H}^{i}(X;\Z_{p}) = 0$ and that $Y$ is a finite-dimensional cohomology 
   $k$-dimensional sphere over the field $\Z_{p}$.
   
   \noindent
   If there exists a $(\Z_{r})^{t}$-equivariant map $f: X {\longrightarrow} Y$, then $\ell < k$.
\end{theorem}

\noindent
\textbf{Box complexes.}
Alon, Frankl, and Lov{\'a}sz~\cite{AFL}, K{\v r}{\' i}{\v z}~\cite{K92}, and Matou{\v s}ek and Ziegler~\cite{MZ} describe
different versions of a box complex to obtain topological lower bounds for a (hyper)graph. We now define a box complex
$\Boxcomplex \s$ associated to an $r$-uniform hypergraph $\s$ which reduces in the graph case to the box complexes $\Boxcomplex G$
described by Matou{\v s}ek and Ziegler~\cite{MZ}. This complex is different from the complexes considered by Alon, Frankl, and 
Lov{\'a}sz~\cite{AFL} and K{\v r}{\'i}{\v z}~\cite{K92}. For an $r$-uniform hypergraph $\s$ without multiplicities we define
\[
\Boxcomplex \s  := \set{ (U_{1},{\ldots} ,U_{r}) \subseteq \vertices \s^{r} }
                       { \begin{matrix}
			    \text{$\bigcup_{i \in [r]} U_{i} \neq \emptyset$ and if all $U_{i} \neq \emptyset$ we have:}\\
			    \text{ $u_{i} \in U_{i}$ ($1 \leq i \leq r$) implies $\{u_{1},{\ldots} ,u_{r}\} \in \edges\s$ }
			 \end{matrix}
                       }.
\]
In case of an $r$-uniform hypergraph with multiplicities we replace the $r$-sets in the definition by $r$-multisubsets.
The box complex of an $r$-uniform hypergraph has a free action by cyclic shift if~$r$ is prime and a fixed-point free action
for arbitrary $r$. 
Since $\kg{r}{s}{\s}$ (considered as a hypergraph with multiplicities) is a subhypergraph of $\KG{r}{s}{\s}$ we have
for prime $r$:
\[
  \zind{r}{\Boxcomplex {\kg{r}{s}{\s}}} \leq \zind{r}{\Boxcomplex {\KG{r}{s}{\s}}}.
\]

\noindent
\textbf{Colour complexes.}\ \
The colour complex $\mathcal C$ of an $r$-uniform hypergraph $\s$ is the simplicial complex $(\sigma^{\chr \s -1})^{*r}_{\Delta(r)}$,
i.e. the vertex set $\vertices {\mathcal C}$ consists of $r$ copies of $[\chr \s]$ and the simplices are ordered $r$-tuples 
$(C_{1}, {\ldots} ,C_{r})$ of subsets of $[\chr \s]$ with $\bigcup_{i \in [n]}C_{i} \neq \emptyset$ and 
$\bigcap_{j=1}^{r}C_{j} = \emptyset$. The colour complex is free if $r$ is prime.
\section{A topological lower bound for the chromatic number}
\label{sec_prime}

\begin{theorem} \label{thm_primepower}
     Let $r=p^{t}$ for a prime $p$ and a positive integer $t$. Consider an $r$-uniform hypergraph~$\s$ with or without multiplicities.
     Suppose that ${\widetilde H}^{i}(\Boxcomplex \s;\Z_{p}) = 0$ for $i \leq \ell$. Then
     \[
       \chr \s \geq \left\lceil \frac{\ell + 2}{r-1} \right\rceil.
     \]
\end{theorem}
\begin{proof}
   A colouring $c: \vertices \s {\longrightarrow} [\chr \s]$ induces a continuous $\Z_{r}$-map 
   \[
     f_{c}: \Boxcomplex \s {\longrightarrow} (\R^{\chr \s})^{r} \setminus \{(v,{\ldots} ,v) \ | \ v\in\R^{\chr \s}\}.
   \]
   Consider the standard basis $e_{1}, {\ldots} , e_{r \cdot \chr \s}$. Map a vertex
   $(\emptyset,{\ldots} ,\emptyset, v, \emptyset,{\ldots} ,\emptyset)$ that has non-empty coordinate $j$ to
   $e_{j \cdot c(v)}$, and extend by linearity. The image of 
   this map is contained in the boundary of the simplex that is spanned by $e_{1}, {\ldots} ,e_{r\cdot \chr \s}$ and
   does not meet $\{(v,{\ldots} ,v) \ | \ v\in\R^{\chr \s}\}$. In particular, a $(\Z_{p})^{t}$-homotopic copy of $\im {f_{c}}$ 
   is contained in a sphere of dimension $((r-1)\cdot \chr \s - 1)$ by normalising each point of $\im {f_{c}}$. The spaces 
   $\Boxcomplex \s$ and $\im {f_{c}}$ (as well its homotopic copy) are fixed-point free $(\Z_{p})^{t}$-spaces, hence we can apply 
   Volovikov's theorem (Theorem~\ref{thm_Volovikov}) to deduce $\ell < (r-1)\cdot \chr \s - 1$.
\end{proof}

\noindent
So far it is not possible to prove this result for arbitrary $r$. In Section~\ref{sec_loops}, we only need the following 
weaker statement. We include its proof since we can avoid Volovikov's theorem.

\begin{theorem} \label{thm_prime_case}
       Let $r$ be a prime and $\s$ be an $r$-uniform hypergraph with or without multiplicities. Then
     \[
       \chr \s \geq \left\lceil \frac{\zind{r}{\Boxcomplex \s} + 1}{r-1} \right\rceil .
     \]
\end{theorem}
\begin{proof}
   We have to show that $\zind{r}{\Boxcomplex \s} \leq (r-1)\cdot \chr \s - 1$. A colouring 
   $c: \vertices \s {\longrightarrow} [\chr \s]$ induces a map $f_{c}: \Boxcomplex \s {\longrightarrow} \mathcal C$ 
   defined on the vertices by
   \[
     (\emptyset,{\ldots} ,\emptyset, v, \emptyset,{\ldots} ,\emptyset)
     \longmapsto
     (\emptyset,{\ldots} ,\emptyset, c(v), \emptyset,{\ldots} ,\emptyset),
   \]
   where $v \in \vertices \s$ and a vertex of $\Boxcomplex \s$ that has non-empty entry in coordinate $i$ is mapped
   to a vertex of $\mathcal C$ that has non-empty entry in coordinate $i$. This map is well-defined since $c$ is a colouring,
   it extends naturally to a simplicial map, and it is $\Z_{r}$-equivariant, no matter whether $\s$ is multiplicity-free or not.
   Hence $\im {f_{c}}$ is a simplicial $\Z_{r}$ subcomplex of $\mathcal C$ and we have
   \[
     \zind{r}{\Boxcomplex \s} \leq \zind{r}{\im {f_{c}}} \leq \Dim {\im {f_{c}}} \leq (r-1)\cdot \chr \s -1
   \]
   since a maximal simplex of $\im {f_{c}}$ contains at most $(r-1)\cdot \chr \s$ many vertices.
\end{proof}

\section{A combinatorial lower bound for the chromatic number of Kneser hypergraphs with multiplicities}
\label{sec_loops}

\begin{theorem}\label{thm_loop_kneser}
  For integers $1 \leq s < r$ and a family $\T$ of subsets of $[n]$, we have
  \[
    \chr {\KG{r}{s}{\T}}
    \geq
    \left\lceil
         \frac{\opcd{r}{s}{\T}}
	      {r-1}
    \right\rceil .
  \]
  If $r$ is prime we have 
  \[
    \chr {\KG{r}{s}{\T}}
    \geq
    \left\lceil
         \frac{\zind{r}{\Boxcomplex {\KG{r}{s}{\T}}}}
	      {r-1}
    \right\rceil
    \geq
    \left\lceil
         \frac{\opcd{r}{s}{\T}}
	      {r-1}
    \right\rceil .
  \]
\end{theorem}

The proof consists of two steps: First one has to prove $\opcd{r}{s}{\T} \leq \zind{r}{\Boxcomplex {\KG{r}{s}{\T}}}$ for
all prime numbers $r$, then the claim follows by induction from Theorem~\ref{thm_prime_case} as shown by Sarkaria~\cite{S90}
or alternatively Ziegler~\cite{Z}. The idea of this induction traces back to Alon, Frankl, and Lov{\' a}sz~\cite{AFL} and was 
also applied by K{\v r}{\' i}{\v z}~\cite{K00} and Matou{\v s}ek~\cite{M02}.

\begin{proof}
   For a subset $U \subseteq [n]$ we define a map $g: 2^{[n]} {\longrightarrow} 2^{\T}$ by 
   $U \longmapsto \set{T\in \T}{T \subseteq U}$. This map is used to define another map 
   $f: \RS{r}{s}{\complex P}\T {\longrightarrow} \sd {\Boxcomplex {\KG{r}{s}{\T}}}$ via
   $(U_{1}, {\ldots} ,U_{r}) \longmapsto (g(U_{1}),{\ldots} , g(U_{r}))$. This map is well-defined since at least one
   $U_{i}$ contains an element of $\T$ and $U^{\prime}_{1}, {\ldots} ,U^{\prime}_{r}$ are $s$-disjoint if 
   $U^{\prime}_{i} \subseteq U_{i}$ and $U_{1}, {\ldots} ,U_{r}$ is $s$-disjoint. The map yields a simplicial map
   because chains of elements of $\RS{r}{s}{P}$ are mapped to chains of simplices of $\Boxcomplex {\KG{r}{S}{\T}}$. Finally,
   the map is $\Z_{r}$-equivariant and surjective, hence
   \[
     \zind{r}{\RS{r}{s}{\complex P}\T} \leq \zind{r}{\im f} 
                                         =  \zind{r}{\sd {\Boxcomplex {\KG{r}{s}{\T}}}}
                                         =  \zind{r}{\Boxcomplex {\KG{r}{s}{\T}}}.
   \]
   To apply Sarkaria's inequality, consider the subcomplex $\complex L$ of $\RS{r}{s}{\complex P}$ that is induced from
   the vertices $\vertices {\RS{r}{s}{\complex P}} \setminus \vertices {\RS{r}{s}{\complex P}\T}$ and use that
   $\RS{r}{s}{\complex P} \subseteq \RS{r}{s}{\complex P}\T * \complex L$. Hence
   \[
     \zind{r}{\RS{r}{s}{\complex P}} \leq \zind{r}{\RS{r}{s}{\complex P}\T * \complex L}
                                     \leq \zind{r}{\RS{r}{s}{\complex P}\T} + \zind{r}{\complex L} + 1.
   \]   
   Since $\zind{r}{\RS{r}{s}{\complex P}} = ns -1$ and since the dimension is an upper bound for the index, 
   we have
   \[
     ns - 1 - \Dim {\complex L} - 1 \leq \zind{r}{\RS{r}{s}{\complex P}\T}.
   \]
   But $\Dim {\complex L} \leq ns - \opcd{r}{s}{\T} - 1$ since every chain of length larger than 
   $ns - \opcd{r}{s}{\T}$ in $\RS{r}{s}{P}$ contains at least one $s$-disjoint $r$-tuple $(U_{1}, {\ldots} , U_{r})$
   that satisfies $\sum_{i=1}{r}|U_{i}| \geq ns - \opcd{r}{s}{\T} + 1$.
\end{proof}

\noindent
We could hide the box complex in the proof and compose 
$f: \RS{r}{s}{\complex P}\T {\longrightarrow} \sd {\Boxcomplex {\KG{r}{s}{\T}}}$ with the ``barycentric subdivision''
$\sd {f_{c}}: \sd{\Boxcomplex \s} {\longrightarrow} \sd {\mathcal C}$ to obtain 
$(r-1)\cdot \chr {\KG{r}{s}{\T}} \geq \opcd{r}{s}{\T}$.

\end{document}